\documentclass[twocolumn,superscriptaddress,showpacs,showkeys]{revtex4}
\usepackage{bm}

\usepackage{amsmath}
\usepackage{amssymb}
\usepackage{amsthm}

\newtheorem{theorem}{Theorem}
\newtheorem{proposition}[theorem]{Proposition}

\DeclareMathOperator{\sgn}{sgn}

\DeclareMathOperator{\curl}{curl}

\begin{document}

\title{Global regularity for a Birkhoff-Rott-$\alpha $ approximation
of the dynamics of vortex sheets of the 2D Euler equations}

% First author block:
\author{Claude Bardos}
\email{bardos@ann.jussieu.fr}
% \homepage{An author's web page; optional}
\affiliation{Universit\'{e} Denis Diderot and Laboratory J.-L.~Lions\\
Universit\'{e} Pierre et Marie Curie, Paris, France }

% For other authors please repeat the author block as needed
\author{Jasmine S.~Linshiz}
\email{jasmine.tal@weizmann.ac.il}
\affiliation{Department of
Computer Science and Applied Mathematics \\Weizmann Institute of
Science \\Rehovot 76100, Israel}

% For other authors please repeat the author block as needed
\author{Edriss S.~Titi}
\email{etiti@math.uci.edu} \email{edriss.titi@weizmann.ac.il}
% Note how REVTeX 4 deals with identical affiliations
\affiliation{Department of Computer Science and Applied
Mathematics \\Weizmann Institute of Science \\Rehovot 76100,
Israel}
\affiliation{Department of Mathematics and\\
Department of Mechanical and Aerospace Engineering \\
University of California \\
Irvine, CA 92697-3875, USA}

\begin{abstract}
We present an $\alpha $-regularization of the Birkhoff-Rott
equation, induced by the two-dimensional Euler-$\alpha$ equations,
for the vortex sheet dynamics. We show that initially smooth
self-avoiding vortex sheet remains smooth for all times under the
$\alpha $-regularized dynamics, provided the initial density of
vorticity is an integrable function over the curve with respect to
the arc-length measure.
\end{abstract}

% Insert suggested PACS numbers (up to 4) in braces.
% The PACS (Physics and Astronomy Classification Scheme)
% can be accessed on the web at http://www.aip.org/pacs/
%47.32.C         Vortex dynamics
%47.20.Ma        Interfacial instabilities (e.g., Rayleigh-Taylor)
%47.20.Ft        Instability of shear flows (e.g., Kelvin-Helmholtz)
%47.15.ki       Inviscid flows with vorticity

% Insert keywords (up to about 4) in braces; optional.

\maketitle

\section{\label{sec:intro}Introduction}
One of the novel approaches for subgrid scale modeling is the
$\alpha $-regularizations of the Navier-Stokes equations (NSE). The
inviscid Euler-$ \alpha $ model was originally introduced in the
Euler-Poincar\'{e} variational framework in \cite{a_HMR98,a_HMR98b}.
In \cite{a_CFHOTW98,a_CFHOTW99, a_CFHOTW99_ChanPipe,a_FHT01,a_FHT02}
the corresponding Navier-Stokes-$ \alpha $ (NS-$\alpha $) [also
known as the viscous Camassa-Holm equations or the
Lagrangian-averaged Navier-Stokes-$\alpha $ (LANS-$\alpha $)] model,
was obtained by introducing the appropriate viscous term into the
Euler-$\alpha $ equations. The extensive research of the $\alpha
$-models (see, e.g., \cite%
{a_HT05,a_CHOT05,a_ILT05,a_VTC05,a_CTV05,a_FHT02,a_FHT01,a_CFHOTW99,a_CFHOTW99_ChanPipe,
a_CFHOTW98,a_MKSM03,a_LL06,a_LL03,a_L06,a_GH03,a_GH06,a_CLT06,a_BFR80,a_HN03,a_CHT04,a_CFR79,a_LT07}%
) stems from the successful comparison of their steady state
solutions to empirical data, for a large range of huge Reynolds
numbers, for turbulent flows in infinite channels and pipes. On the
other hand, the $\alpha $-models can also be viewed as numerical
regularizations of the original, Euler or Navier-Stokes, systems.
The main practical question arising is that of the applicability of
these regularizations to the correct predictions of the underlying
flow phenomena.

In this paper we present some results concerning the $\alpha
$-regularization of the two-dimensional (2D) Euler equations in the
context of  vortex sheet dynamics. A vortex sheet is a surface of
codimension one (a curve in the plane) in inviscid incompressible
flow, across which the tangential component of the velocity has a
jump discontinuity, while the normal component is continuous. The
evolution of the vortex sheet can be described by the Birkhoff-Rott
(BR) equation \cite{a_B62,a_R56,b_S92}. This is a nonlinear singular
integro-differential equation, which can be obtained formally from
the Euler equations assuming that the evolution of a vortex sheet
retains a curve-like structure. However, the initial data problem
for the BR equation is ill-posed due to the Kelvin-Helmholtz
instability \cite{a_B62, a_SB79}. Numerous results show that an
initially real analytic vortex sheet can develop a finite time
singularity in its curvature. This singularity formation was studied
with asymptotic techniques in \cite{a_M79,a_CBT00} and numerically
in \cite{a_MBO82,a_K86a,a_CBT00}. Specific examples of solutions
were constructed in \cite{a_DR88,a_CO89}, where the development, in
a finite time, of curvature singularity from initially analytic data
was rigorously proved.

The problem of the evolution of a vortex sheet can also be
approached, in the general framework of weak solutions (in the
distributional sense) of the Euler equations, as a problem of
evolution of the vorticity, which is concentrated as a measure along
a surface of codimension one. The general problem of existence for
mixed-sign vortex sheet initial data remains an open question.
However, in 1991, Delort \cite{a_D91} proved a global in time
existence of weak solutions of the 2D incompressible Euler equation
for the vortex sheet initial data with initial vorticity being a
Radon measure of a distinguished sign, see also \cite
{a_EM94,a_M93,a_LX95, a_S95, a_S96, b_MB02}. This result was later
obtained as an inviscid limit of the Navier-Stokes regularizations
of the Euler equations \cite{a_M93, a_S95}, and as a limit of vortex
methods \cite{a_LX95, a_S96}. The Delort's result was also extended
to the case of mirror-symmetric flows with distinguished sign
vorticity on each side of the mirror \cite{a_LfNlX01}. However, the
problem of uniqueness of a weak solution with a fixed sign vortex
sheet initial data is still unanswered, numerical evidences of
non-uniqueness can be found, e.g., in \cite{a_P89, a_FLLZ06}.
Furthermore, the structure of weak solutions given by Delort's
theorem is not known, while the Birkhoff-Rott equations assume
\textit{a priori} that a vortex sheet remains a curve at a later
time. A proposed criterion for the equivalence of a weak solution of
the 2D Euler equations with vorticity being a Radon measure
supported on a curve, and a weak solution of the Birkhoff-Rott
equation can be found in \cite{a_FLS06}. Also, another definition of
weak solutions of Birkhoff-Rott equation has been proposed in
\cite{a_W02,a_W06}. For a recent survey of the subject, see
\cite{a_BT07}.

The question of global existence of weak solutions for the
three-dimensional Euler-$\alpha $ equations is still an open
problem. On the other hand, the 2D Euler-$\alpha $ equations were
studied in \cite{a_OS01}, where it has been shown that there exists
a unique global weak solution to the Euler-$ \alpha $ equations with
initial vorticity in the space of Radon measures on $
{\mathbb{R}}^{2}$, with a unique Lagrangian flow map describing the
evolution of particles. In particular, it follows that the
vorticity, initially supported on a curve, remains supported on a
curve for all times.

We present in this paper an analytical study of the $\alpha
$-analogue of the Birkhoff-Rott equation, the
Birkhoff-Rott-$\alpha$ (BR-$\alpha $) model, which is induced by
the 2D Euler-$\alpha$ equations. The BR-$\alpha $ model was
implemented computationally in \cite{a_HNP06}, where a numerical
comparison between the BR-$\alpha $ regularization and the
existing regularizing methods, such as a vortex blob model
\cite{a_CB73, a_K86b, a_CK00,a_LX95,a_BP06} has been performed. We
remark that, unlike the vortex blob methods that regularize the
singular kernel in the Birkhoff-Rott equation, the $\alpha$-model
regularizes instead the Euler equations themselves to obtain a
smoother kernel.

We report in Section~\ref{sec:GlobalReg} our main result, which
states that the initially smooth self-avoiding 2D vortex sheet,
evolving under the BR-$\alpha $ equation, remains smooth for all
times. In this short communication we only report the results and
sketch some of their proofs, the full details will be reported in a
forthcoming paper. In Section~\ref{sec:BR_alpha} we describe the
BR-$\alpha$ equation. Section~\ref{sec:LinStab} studies the linear
stability of a flat vortex sheet with uniform vorticity density for
the 2D BR-$\alpha $ model. The linear stability analysis shows that
the BR-$\alpha $ regularization controls the growth of high wave
number perturbations, which is the reason for the well-posedness.
This is unlike the case for the original BR problem that exhibits
the Kelvin-Helmholtz instability, the main mechanism for its
ill-posedness.

\section{\label{sec:BR_alpha}Birkhoff-Rott-$\protect\alpha $ equation}

The incompressible Euler equations in $\mathbb{R}^{2}$ in the
vorticity form are given by
\begin{equation}
\begin{split}
& \frac{\partial q}{\partial t}+\left( v\cdot \nabla \right) q=0, \\
& v=K\ast q, \\
& q(x,0)=q^{in}(x),
\end{split}
\label{grp:EulerEqVortForm}
\end{equation}
where $K\left( x\right) =\frac{1}{2\pi }\nabla ^{\perp }\log
\left\vert x\right\vert $, $v$ is the fluid velocity field,
\mbox{$q=\curl v$} is the vorticity, and $q^{in}$ is the given
initial vorticity.

The 2D Euler-$\alpha $ model \cite%
{a_CFHOTW99_ChanPipe,a_HMR98,a_HMR98b,a_H02_pA,a_MS03,a_C01} is an
inviscid regularization of the Euler equations, such that the
vorticity is governed by the system
\begin{equation}
\begin{split}
& \frac{\partial q}{\partial t}+\left( u\cdot \nabla \right) q=0, \\
& u=K^{\alpha }\ast q, \\
& q(x,0)=q^{in}(x).
\end{split}
\label{grp:Euler_alpha_vortForm}
\end{equation}
Here $u$ represents the ``filtered'' fluid velocity, and $\alpha
>0$ is a length scale parameter, which represents the width of the
filter. At the limit $\alpha =0$, we formally obtain the Euler
equations \eqref{grp:EulerEqVortForm}. The smoothed kernel is
$K^{\alpha }=G^{\alpha }\ast K$, where $G^{\alpha }$ is the Green
function associated with the Helmholtz operator $\left( I-\alpha
^{2}\Delta \right) $, given by
\begin{equation}
G^{\alpha }\left( x\right) =\frac{1}{\alpha ^{2}}G\left( \frac{x}{\alpha }%
\right) =-\frac{1}{\alpha ^{2}}\frac{1}{2\pi }K_{0}\left(
\frac{\left\vert x\right\vert }{\alpha }\right),
\label{eq:GreenFunc_Helmholtz_2D}
\end{equation}
here $x=\left( x_{1},x_{2}\right) \in \mathbb{R}^{2}$ and $K_{0}$ is
a modified Bessel function of the second kind \cite{b_W44}.

Let $\mathcal{M}({\mathbb{R}}^{2})$ denote the space of Radon
measures on ${ \mathbb{R}}^{2}$; $\mathcal{G}$ denote the group of
all homeomorphism of ${ \mathbb{R}}^{2}$, which preserve the
Lebesgue measure; and $\eta=\eta (\cdot ,t)$ denote the Lagrangian
flow map induced by \eqref{grp:Euler_alpha_vortForm} and obeying the
equation $\partial _{t}\eta (x,t)=u(\eta (x,t),t)$, \, $\eta \left(
x,0\right) =x$.

Oliver and Shkoller \cite{a_OS01} showed global well-posedness of
the Euler-$ \alpha $ equations \eqref{grp:Euler_alpha_vortForm} with
initial vorticity in $\mathcal{M}({\mathbb{R}}^{2})$ (which includes
point-vortex data).

\begin{theorem}
\label{thm:OS01}\emph{(Oliver and Shkoller \cite{a_OS01})} For
initial data $q^{in}\in \mathcal{M}({\mathbb{R}}^{2})$, there exists
a unique global weak solution (in the sense of distribution) to
\eqref{grp:Euler_alpha_vortForm} with
\begin{equation*}
\eta \in C^{1}\left( {\mathbb{R}};\mathcal{G}\right) ,\,u\in C\left(
{ \mathbb{R}};C\left( {\mathbb{R}}^{2}\right) \right),\,
% \text{and }
q\in C\left( {\mathbb{R}};\mathcal{M}({\mathbb{R}}^{2})\right).
\end{equation*}%
\end{theorem}

The Birkhoff-Rott-$\alpha $ equation, based on the Euler-$ \alpha $
equations is derived similarly to the derivation of the original
Birkhoff-Rott equation. Detailed descriptions of the Birkhoff-Rott
equation as a model for the evolution of the vortex sheet can be
found, e.g., in \cite{b_S92,b_MP94,b_MB02}. We remark, that while
the BR equations assume \textit{a priori} that a vortex sheet
remains a curve at a later time, in the 2D Euler-$\alpha$ case, if
the vorticity is initially supported on a curve, then due to the
existence of the unique Lagrangian flow map given by Theorem
\ref{thm:OS01}, it remains supported on a curve for all times. Hence
the BR-$\alpha$ equation gives an equivalent description of the
vortex sheet evolution, as the 2D Euler-$\alpha$ equations. It is
described in the following proposition.

\begin{proposition}
Let $q$ be the solution of \eqref{grp:Euler_alpha_vortForm} in the
sense of the Theorem \ref{thm:OS01}. Assume, furthermore, that $q$
has the density $\gamma (\sigma ,t)$ supported on the sheet (curve)
$\Sigma (t)=\left\{ x=x(\sigma ,t)\in \mathbb{R}^{2}|\sigma
_{0}\left( t\right) \leq \sigma \leq \sigma _{1}\left( t\right)
\right\} $, that is, the vorticity $q(x,t)$ satisfies
\begin{equation*}
\int_{{\mathbb{R}}^{2}}\varphi (x)dq(x,t)=\int_{\sigma _{0}\left(
t\right) }^{\sigma _{1}\left( t\right) }\varphi \left( x(\sigma
,t)\right) \gamma (\sigma ,t)|x_{\sigma }\left(\sigma ,t
\right)|d\sigma ,
\end{equation*}
for every $\varphi \in C_{0}^{\infty }\left(
{\mathbb{R}}^{2}\right)$, $\gamma \left( \cdot ,t\right) \in
L^{1}(\left\vert x_{\sigma }\right\vert d\sigma )$. Then this sheet
evolves according to the equation
\begin{align*}
&\frac{\partial }{\partial t}x\left( \sigma ,t\right)= \\
&=\int_{\sigma _{0}\left( t\right) }^{\sigma _{1}\left( t\right)
}K^{\alpha }\left( x\left( \sigma ,t\right) -x\left( \sigma ^{\prime
},t\right) \right) \gamma \left( \sigma ^{\prime },t\right)
\left\vert x_{\sigma }\left(\sigma^{\prime } ,t \right)\right\vert
d\sigma ^{\prime }.
\end{align*}
Additionally, if $\Gamma \left( \sigma ,t\right) =\int_{\sigma
^{\ast }}^{\sigma }\gamma \left( \sigma ^{\prime },t\right)
\left\vert x_{\sigma }\left( \sigma ^{\prime },t\right) \right\vert
d\sigma ^{\prime }$, where $ x\left( \sigma ^{\ast },t\right) $ is
some fixed reference point on $\Sigma (t)$, defines a strictly
increasing function of $\sigma $ (e.g., as in the case of positive
vorticity), then the evolution equation is given by the
Birkhoff-Rott-$\alpha $ (BR-$\alpha $) equation
\begin{equation}  \label{eq:BR_alpha}
\frac{\partial }{\partial t}x\left( \Gamma ,t\right) =\int_{\Gamma
_{0}}^{\Gamma _{1}}K^{\alpha }\left( x\left( \Gamma ,t\right)
-x\left( \Gamma ^{\prime },t\right) \right) d\Gamma ^{\prime }
\end{equation}
with $\gamma =1/|x_{\Gamma }|$ being the vorticity density along the
sheet and $-\infty <\Gamma _{0}<\Gamma _{1}<\infty $.
\end{proposition}

Here $\sigma _{0},\sigma _{1}$ can represent either a finite length
curve, or an infinite one. Notice that
\begin{equation*}
K^{\alpha }\left( x\right) =\nabla ^{\perp }\Psi ^{\alpha }\left(
\left\vert x\right\vert \right) =\frac{x^{\perp }}{\left\vert
x\right\vert }D\Psi ^{\alpha }\left( \left\vert x\right\vert
\right),
\end{equation*}
where
\begin{equation*}
\Psi ^{\alpha }\left( r\right) =\frac{1}{2\pi }%
\left[ K_{0}\left( \frac{r }{\alpha }\right) +\log r \right]
\end{equation*}
and
\begin{equation*}
D\Psi ^{\alpha }(r)=\frac{d\Psi ^{\alpha }}{dr}(r) =\frac{1}{2\pi }%
\left[ -\frac{1}{\alpha }K_{1}\left( \frac{ r }{\alpha }%
\right) +\frac{1}{ r }\right].
\end{equation*}
$K_{0}$ and $K_{1}$ denote modified Bessel functions of the second
kind of orders zero and one, respectively. For details on Bessel
functions, see, e.g., \cite{b_W44}. We remark that the smoothed
kernel $K^{\alpha }\left( x\right) $ is a
bounded continuous function, that for $\frac{\left\vert x\right\vert }{%
\alpha }\rightarrow 0$ behaves as $K^{\alpha }\left( x\right) =-\frac{1}{%
4\pi }\frac{1}{\alpha ^{2}}x^{\perp }\log \frac{\left\vert x\right\vert }{%
\alpha }+O\left( \frac{|x|}{\alpha ^{2}}\right) $. That is, it is
non-singular kernel. Since $\gamma \left( \cdot ,t\right) \in
L^{1}(\left\vert x_{\sigma }\right\vert d\sigma )$ we can show the
integrability of the relevant terms, even though $\left\vert
K^{\alpha }\left( x\right) \right\vert $ is decaying like
$\left\vert x\right\vert ^{-1}$ at infinity.

\section{\label{sec:LinStab}Linear stability of a flat vortex sheet with uniform vorticity
density for 2D BR-$\protect\alpha $ model}

The initial data problem for the BR equation is highly unstable due
to an ill-posed response to small perturbations called
Kelvin-Helmholtz instability \cite{a_B62, a_SB79}. The linear
stability analysis of the BR-$\alpha$ equation shows that the
ill-posedness of the original problem is mollified, and the
Kelvin-Helmholtz instability of the original system now disappears.

When the vortex sheet can be parameterized as a graph of a function
in the form $x_{2}=x_{2}\left( x_{1},t\right) $ the BR-$\alpha $
system \eqref{eq:BR_alpha} takes the form
\begin{align}
\frac{\partial x_{2}}{\partial t}& =-\frac{\partial x_{2}}{\partial x_{1}}%
u_{1}+u_{2},  \label{eq:2dSystem} \\
\frac{\partial \gamma }{\partial t}& =-\frac{\partial }{\partial x_{1}}%
\left( \gamma u_{1}\right) ,  \notag
\end{align}
with velocity $u =\left( u_{1},u_{2}\right) ^{t}$ given by
\begin{align*}
u\left( x_{1},t\right)
=\mathrm{p.v.}\int_{\mathbb{R}}K^{\alpha}\left( x\left(
x_{1},t\right) -x\left( x_{1}^{\prime },t\right) \right) \gamma
\left( x_{1}^{\prime },t\right) dx_{1}^{\prime },
\end{align*}
where  $x\left( x_{1},t\right) =\left( x_{1},x_{2}\left(
x_{1},t\right) \right) ^{t}$. The flat sheet \mbox{$x_{2}^{0}\equiv
0$} with uniformly concentrated intensity $ \gamma _{0}$ is
stationary solution of \eqref{eq:2dSystem}. By linearization about
the flat sheet we obtain the following linear system
\begin{align*}
& \frac{\partial \tilde{x}_{2}}{\partial t}=\tilde{u}_{2}, \\
& \frac{\partial \tilde{\gamma}}{\partial t}=-\gamma
_{0}\frac{\partial \tilde{u}_{1}}{\partial x_{1}},
\end{align*}
where
\begin{align*}
\tilde{u}_{1}\left( x_{1},t\right) & =-\gamma _{0}\left(
\sgn\left( x_{1}\right) D\Psi ^{\alpha }\left( \left\vert
x_{1}\right\vert \right)
\right) \ast \frac{\partial \tilde{x}_{2}}{\partial x_{1}}, \\
\tilde{u}_{2}\left( x_{1},t\right) & =\left( \sgn\left( x_{1}\right)
D\Psi ^{\alpha }\left( \left\vert x_{1}\right\vert \right) \right)
\ast \tilde{ \gamma},
\end{align*}
and $\left( \tilde{x}_{2},\tilde{\gamma}\right) $ is a small
perturbation about the flat sheet.

Consequently, the equation for the Fourier modes is given by
\begin{equation}
\frac{d}{dt}
\begin{pmatrix}
\widehat{\tilde{x}_{2}} \\
\widehat{\tilde{\gamma}}
\end{pmatrix}
=
\begin{pmatrix}
0 & \frac{i}{2}\sgn(k)d(k) \\
-i\frac{\gamma _{0}^{2}}{2}k^{2}\sgn(k)d(k) & 0
\end{pmatrix}
\begin{pmatrix}
\widehat{\tilde{x}_{2}} \\
\widehat{\tilde{\gamma}}
\end{pmatrix}
,  \label{eq:FourierModesSystem}
\end{equation}
where
\begin{equation*}
d(k)= \left( 1+\frac{1}{\alpha ^{2}k^{2}}\right) ^{-1/2}-1.
\end{equation*}
Observe that in order to calculate the Fourier transform
\begin{equation*}
\mathcal{F}\left( \sgn\left( x_{1}\right) D\Psi ^{\alpha }\left(
\left\vert x_{1}\right\vert \right) \right) \left( k\right)
=\frac{i}{2}\sgn(k)d(k) ,
\end{equation*}
we used the integral representation of the modified Bessel function
of the second kind $K_{1}\left( x_{1}\right) =x_{1}\int_{1}^{\infty
}e^{-x_{1}t}\left( t^{2}-1\right) ^{1/2}dt$, (see, e.g.,
\cite{b_W44}). The eigenvalues of the coefficient matrix, given in
\eqref{eq:FourierModesSystem}, are
\begin{equation}\label{eq:eigenval}
\lambda(k) =\pm \frac{1}{2}\left\vert \gamma _{0}\right\vert
\left\vert k\right\vert \left( 1-\left( 1+\frac{1}{\alpha
^{2}k^{2}}\right) ^{-1/2}\right) .
\end{equation}
To conclude, the $\alpha$-regularization mollifies the
Kelvin-Helmholtz instability as follows: we have an algebraic decay
of the eigenvalues to zero of order $ \frac{1}{\alpha ^{2}\left\vert
k\right\vert }$, as $k\rightarrow \infty $ ($ \alpha $ fixed).
While, for $\alpha $ $\rightarrow 0$, for fixed $k$, we recover the
eigenvalues of the original BR equations $\pm \frac{1}{2} \left\vert
\gamma _{0}\right\vert \left\vert k\right\vert $ (see, e.g., \cite
{a_SSBF81}).

For the sake of comparison, we observe that for the vortex blob
regularization of Krasny \cite{a_K86a}, where the singular BR
kernel, $K(x)$, was replaced with the smoothed kernel
\begin{equation*}
K_{\delta }\left( x\right) =K\left( x\right) \frac{\left\vert
x\right\vert ^{2}}{\left\vert x\right\vert ^{2}+\delta
^{2}}=\frac{1}{2\pi }\frac{ x^{\perp }}{\left\vert x\right\vert
^{2}+\delta ^{2}},
\end{equation*}
the eigenvalues are
\begin{equation*}
\lambda(k) =\pm \frac{1}{2}e^{-\delta k}\left\vert \gamma
_{0}\right\vert \left\vert k\right\vert
\end{equation*}
with an exponential decay to zero, as $k\rightarrow \infty $
($\delta>0 $ is fixed). As $\delta $ $\rightarrow 0$, for fixed $k$,
one recovers again the eigenvalues of the original BR equations.

The behavior of the eigenvalues of the linearized system
\eqref{eq:FourierModesSystem} indicates that high wave number
perturbations grow exponentially in time with a rate that decays to
zero, as $k\rightarrow \infty $, which is the reason for
well-posedness of the $\alpha$-regularized model. This is unlike the
original BR problem that exhibits the Kelvin-Helmholtz instability.
It is worth mentioning that the $\alpha$-regularization is ``closer"
to the original system than the vortex-blob method at the high wave
numbers, due to the algebraic decay instead of exponential one in
the vortex blob method. This result was also evaluated
computationally in \cite{a_HNP06}.

\section{\label{sec:GlobalReg}Global regularity for BR-$\protect\alpha $ equation}

In this section we present the global existence and uniqueness of
solutions of the BR-$\alpha$ equation \eqref{eq:BR_alpha} in the
appropriate space of functions. We show that initially smooth
solutions of \eqref{eq:BR_alpha} remain smooth for all times.

Let us first describe the H\"{o}lder space $C^{n,\beta }\left(
\Sigma \subset \mathbb{R};\mathbb{R}^{2}\right) $, $0<\beta \leq 1$,
which is the space of functions $x:\Sigma \subset
\mathbb{R}\rightarrow \mathbb{R}^{2}$, with finite norm
\begin{equation*}
\left\Vert x\right\Vert _{C^{n,\beta }\left( \Sigma \right) }=\sum_{
k =0}^{n} \left\vert \frac{d^k}{d{\Gamma}^k} x\right\vert
_{C^{0}\left( \Sigma \right) }+\left\vert \frac{d^n}{d{\Gamma}^n}
x\right\vert _{\beta \left( \Sigma \right) },
\end{equation*}
where
\begin{equation*}
\left\vert x\right\vert _{C^{0}\left( \Sigma \right)
}=\sup_{\Gamma \in \Sigma }\left\vert x\left( \Gamma \right)
\right\vert
\end{equation*}
and $\left\vert \cdot \right\vert _{\beta }$ is the H\"{o}lder
semi-norm
\begin{equation*}
\left\vert x\right\vert _{\beta \left( \Sigma \right) }=\sup_{\substack{ %
\Gamma ,\Gamma ^{\prime }\in \Sigma  \\ \Gamma \neq \Gamma ^{\prime }}}\frac{%
\left\vert x\left( \Gamma \right) -x\left( \Gamma ^{\prime
}\right) \right\vert }{\left\vert \Gamma -\Gamma ^{\prime
}\right\vert ^{\beta }}.
\end{equation*}
We also use the notation
\begin{equation*}
\left\vert x\right\vert _{\ast }=\inf_{\substack{ \Gamma ,\Gamma
^{\prime }\in \Sigma  \\ \Gamma \neq \Gamma ^{\prime
}}}\frac{\left\vert x\left( \Gamma \right) -x\left( \Gamma
^{\prime }\right) \right\vert }{\left\vert \Gamma -\Gamma ^{\prime
}\right\vert }.
\end{equation*}
Next we state our main result.
\begin{theorem}\label{thm:GlobalEx}
Let $n\geq 1$, $0<\beta <1$, $x\left( \Gamma ,0\right) =x_{0}\left(
\Gamma \right) \in C^{n,\beta }\left( \Gamma _{0},\Gamma
_{1}\right)\cap \left\{ \left\vert x\right\vert _{\ast }>0\right\}
$, then for any $T>0$ there is a unique solution \mbox{$x\in
C^{1}\left([-T,T];C^{n,\beta }\left( \Gamma _{0},\Gamma
_{1}\right)\cap \left\{ \left\vert x\right\vert _{\ast
}>0\right\}\right) $} of \eqref{eq:BR_alpha}.

In particular, if \mbox{$x_{0} \in C^{\infty}\left( \Gamma
_{0},\Gamma _{1}\right)\cap \left\{ \left\vert x\right\vert _{\ast
}>0\right\} $} then \mbox{$x\in C^{1}\left([-T,T];C^{\infty}\left(
\Gamma _{0},\Gamma _{1}\right)\cap \left\{ \left\vert x\right\vert
_{\ast }>0\right\}\right) $}.
\end{theorem}
We remark that, although the kernel $K^{\alpha}$ is a continuous
bounded function, its derivatives are unbounded near the origin, and
the condition $|x|_{\ast }>0$, which generally means self-avoiding
curves, allows us to show the integrability of the relevant terms.
Furthermore, it is also worth mentioning that $|x|_{\ast }$ being
bounded away from zero is similar to the chord arc hypothesis
\cite{a_D82}, used later in \cite{a_W02,a_W06} .

Now we sketch the main steps involved in the proof of Theorem
\ref{thm:GlobalEx}. First, we apply the Contraction Mapping
Principle to the BR-$\alpha $ equation \eqref{eq:BR_alpha} to prove
the short time existence and uniqueness of solutions in the
appropriate space of functions. We show that initially $C^{1,\beta
}$ smooth solutions of \eqref{eq:BR_alpha} remain $C^{1,\beta }$
smooth for a finite short time. Next, we derive an \textit{a priori}
bound for the controlling quantity for continuing the solution for
all time. Then we extend the result for higher derivatives. The full
details will be reported in a forthcoming paper.

\begin{proof}[Sketch of the proof]
We consider the BR-$\alpha $ equation as an evolution functional
equation in the Banach space $C^{n,\beta }$
\begin{equation}
\begin{split}
& \frac{\partial x}{\partial t}\left( \Gamma ,t\right)
=\int_{\Gamma _{0}}^{\Gamma _{1}}K^{\alpha }\left( x\left( \Gamma
,t\right) -x\left(
\Gamma ^{\prime },t\right) \right) d\Gamma ^{\prime }, \\
& x\left( \Gamma ,0\right) =x_{0}\left( \Gamma \right) \in
C^{n,\beta }\cap \left\{ \left\vert x\right\vert _{\ast
}>0\right\}
\end{split}
\label{grp:BR_alpha_onBanach}
\end{equation}
with $\gamma =1/|x_{\Gamma }|$ being the vorticity density along
the sheet. Notice that the initial density is well defined for the
subset $\left\{ \left\vert x\right\vert _{\ast }>0\right\} $.

\bigskip \noindent \textit{Step 1.} %
We show the local existence and uniqueness of solutions. To apply
the Contraction Mapping Principle to the BR-$\alpha $ equation
\eqref{grp:BR_alpha_onBanach} we first prove the following
proposition
\begin{proposition}
Let $1<M<\infty $, $-\infty <\Gamma _{0}<\Gamma _{1}$, and let
$S^{M}$ be the set
\begin{equation*}
\left\{ \Gamma \mapsto x\left( \Gamma \right) \in C^{1,\beta }\left(
\Gamma _{0},\Gamma _{1}\right) ,\left\vert x_{\Gamma }\right\vert
_{C^{0}}<M,\left\vert x\right\vert _{\ast }>\frac{1}{M}\right\} .
\end{equation*}
Then the mapping $x\left( \Gamma ,t\right) \mapsto$
\begin{equation*}
u\left( x\left( \Gamma ,t\right) ,t\right) =\int_{\Gamma
_{0}}^{\Gamma _{1}}K^{\alpha }\left( x\left( \Gamma ,t\right)
-x\left( \Gamma ^{\prime },t\right) \right) d\Gamma ^{\prime }
\end{equation*}
defines a locally Lipschitz continuous map from $S^{M}$ into
$C^{1,\beta }$.
\end{proposition}

This implies the local existence and uniqueness of solutions:

\begin{proposition}
Given $x_{0}\left( \Gamma \right) \in C^{1,\beta }\left( \Gamma
_{0},\Gamma _{1}\right) \cap \left\{ \left\vert x\right\vert _{\ast
}>0\right\} $, there
exists $1<M<\infty $ and a time $T(M)$, such that the system %
\eqref{grp:BR_alpha_onBanach} has a unique local solution $x\in
C^{1}((-T(M),T(M));S^{M})$.
\end{proposition}

\bigskip \noindent \textit{Step 2.} %
The obtained local solutions can be continued in time provided that
we have global, in time, bounds on $\frac{1}{\left\vert x\left(
\cdot ,t\right) \right\vert _{\ast }}$ and $\left\vert x_{\Gamma
}\left( \cdot ,t\right) \right\vert _{\beta }$. To control these
quantities we need to bound \mbox{$\int_{0}^{T}\left\Vert \nabla
_{x}u\left( x(\cdot,t) ,t\right) \right\Vert _{L^{\infty }
\left(\Gamma_0,\Gamma_1\right)}dt$}. We sketch the proof of this
bound.
%\begin{proposition}
%\begin{equation*}
%\left\vert \nabla _{x}u\left( x,t\right) \right\vert \leq C\left( \frac{1}{%
%\alpha },t,\left\Vert q^{in}\right\Vert _{\mathcal{M}\left(
%\mathbb{R} ^{2}\right) },\left\vert x_{0}\right\vert _{\ast }\right)
%\end{equation*}
%\end{proposition}
We write $\nabla _{x}u\left( x(\Gamma,t),t\right)$ as
\begin{align*}
\nabla _{x}u\left( x(\Gamma,t),t\right)  &=\int_{\Gamma
_{0}}^{\Gamma _{1}} \nabla_{x} K^{\alpha}\left( x\left( \Gamma
,t\right) -x\left( \Gamma ^{\prime
},t\right) \right)  d\Gamma ^{\prime } \\
& =\int_{E_{\varepsilon }}+\int_{\left( \Gamma _{0},\Gamma
_{1}\right) \backslash E_{\varepsilon }}=I_{1}+I_{2},
\end{align*}
where
\begin{equation*}
 E_{\varepsilon }=\left\{ \Gamma ^{\prime }\in \left( \Gamma
_{0},\Gamma _{1}\right) :\frac{\left\vert x\left( \Gamma ,t\right)
-x\left( \Gamma ^{\prime },t\right) \right\vert }{\alpha
}<\varepsilon \right\},
\end{equation*}
for a fixed small $0<\varepsilon <1$, to be further refined later.
Let $\eta $ denote the unique Lagrangian flow map given by Theorem
\ref{thm:OS01}. Denote the distance between two points $\eta(x,t)$
and $\eta(x',t)$ by $r\left( t\right) =\left\vert \eta \left(
x,t\right) -\eta \left( x^{\prime },t\right) \right\vert $, where
$r\left( 0\right) =\left\vert x-x^{\prime }\right\vert$.

Then, using the estimate (2.14) of \cite{a_OS01}, we have
\begin{align*}
\left\vert \frac{d}{dt}r\left( t\right) \right\vert & \leq
\int_{\mathbb{R} ^{2}}\left\vert K^{\alpha }\left( x,y\right)
-K^{\alpha }\left( x^{\prime
},y\right) \right\vert \left \vert q\left( y,t\right)  \right\vert dy \\
& \leq C\frac{1}{\alpha }\varphi \left( \frac{r\left( t\right)
}{\alpha }
\right) \left\Vert q\right\Vert _{M\left( \mathbb{R}^{2}\right) }\\
&=C\frac{1}{\alpha }\varphi \left( \frac{r\left( t\right) }{\alpha }
\right) \left\Vert q^{in}\right\Vert _{M\left( \mathbb{R}^{2}\right)
},
\end{align*}
where
\begin{equation*}
\varphi \left( r\right) =\left\{
\begin{array}{ll}
0 , & r=0, \\
r\left( 1-\log r\right) , & 0< r<1, \\
1, & r\geq 1.
\end{array}
\right.
\end{equation*}
By comparison with the solution of the differential equation
\begin{equation*}
\frac{d}{dt}r\left( t\right) =-C\frac{1}{\alpha }\varphi \left(
\frac{ r\left( t\right) }{\alpha }\right) \left\Vert
q^{in}\right\Vert _{\mathcal{M} \left( \mathbb{R}^{2}\right) },
\end{equation*}
we can choose $\varepsilon =\varepsilon \left( t,\frac{1}{\alpha }
,\left\Vert q^{in}\right\Vert _{\mathcal{M}\left(
\mathbb{R}^{2}\right) }\right) $ small enough, such that, for
$\frac{\left\vert x\left( \Gamma ,t\right) -x\left( \Gamma ^{\prime
},t\right) \right\vert }{\alpha } <\varepsilon $,
\begin{equation}\label{eq:trajBelowBound}
\begin{split}
& \frac{\left\vert x\left( \Gamma ,t\right) -x\left( \Gamma ^{\prime
},t\right) \right\vert }{\alpha }\geq  \\
& \qquad \geq \left( \frac{\left\vert x\left( \Gamma ,0\right)
-x\left( \Gamma ^{\prime },0\right) \right\vert }{\alpha }\right)
^{e^{tC_{1}}}e^{1-e^{tC_{1}}},
\end{split}
\end{equation}
where $C_{1}=\frac{C}{\alpha ^{2}}\left\Vert q^{in}\right\Vert
_{\mathcal{M} \left( \mathbb{R}^{2}\right) }$. Now, using also that
$|x_0|_{\ast }$ is bounded away from zero, we can bound
$\frac{\left\vert x\left( \Gamma ,t\right) -x\left( \Gamma ^{\prime
},t\right) \right\vert }{\alpha }$ from below, which in turn implies
the bound
\begin{equation*}
I_{1}\leq C\left( t,\frac{1}{\alpha }%
,\left\Vert q^{in}\right\Vert _{\mathcal{M}\left(
\mathbb{R}^{2}\right) },\left\vert x_{0}\right\vert _{\ast }\right)
.
\end{equation*}
While to bound $I_{2}$, we use the boundness of $\left\vert \nabla
_{x}K^{\alpha}\left( x\left( \Gamma ,t\right) ,x\left( \Gamma
^{\prime },t\right) \right) \right\vert$ in $
\{ \Gamma ^{\prime }\in \left( \Gamma _{0},\Gamma _{1}\right) :%
\frac{\left\vert x\left( \Gamma ,t\right) -x\left( \Gamma ^{\prime
},t\right) \right\vert }{\alpha }\geq \varepsilon \}  $. Hence
\begin{align}\label{eq:Du_bound}
\begin{split}
\int_{0}^{T}&\left\Vert \nabla _{x}u\left( x(\cdot,t) ,t\right)
\right\Vert _{L^{\infty } \left(\Gamma_0,\Gamma_1\right)}dt \leq \\
& \qquad \leq C\left( \frac{1}{ \alpha },T,\left\Vert
q^{in}\right\Vert _{\mathcal{M}\left( \mathbb{R} ^{2}\right)
},\left\vert x_{0}\right\vert _{\ast }\right).
\end{split}
\end{align}

Now, by Gr\"{o}nwall inequality the bound \eqref{eq:Du_bound}
provides bounds on $\frac{1}{\left\vert x\left( \cdot ,t\right)
\right\vert _{\ast }}$ and $ \left\vert x_{\Gamma }\left( \cdot
,t\right) \right\vert _{C^{0}}$ on $\left[ 0,T\right] $. The bound
on $\left\vert x_{\Gamma }\left( \cdot ,t\right) \right\vert _{\beta
}$ on $\left[ 0,T\right] $ is a consequence of
\begin{align*}
&\frac{d}{dt}x_{\Gamma }\left( \Gamma ,t\right) =\nabla _{x}u\left(
x\left( \Gamma ,t\right) ,t\right) \cdot x_{\Gamma }\left( \Gamma
,t\right) ,\\
&\left\vert \nabla _{x}u\left(  x\left( \cdot ,t\right) ,t\right)
\right\vert _{\beta }\leq C\left(\frac{1}{\alpha} , \left\vert
x_{\Gamma }\right\vert _{L^{\infty }},\left\vert x\right\vert _{\ast
},\Gamma _{1}-\Gamma _{0}\right) ,
\end{align*}
\eqref{eq:Du_bound}  and the Gr\"{o}nwall inequality.

This yields global in time existence and uniqueness of $C^{1,\beta
}$ solutions of \eqref{grp:BR_alpha_onBanach}.

\bigskip \noindent \textit{Step 3.} %
To provide an a priori bound for higher derivatives in terms of
lower ones, we show that for $x\in S^{M}\cap C^{n,\beta }\left(
\Gamma _{0},\Gamma _{1}\right) $,
\begin{equation*}
\left\vert u\left( x\left( \cdot ,t\right) ,t\right) \right\vert
_{n,\beta }\leq C\left( \frac{1}{\alpha} ,M,\left\vert x\left( \cdot
,t\right) \right\vert _{n-1,\beta }\right) \left\vert x\left( \cdot
,t\right) \right\vert _{n,\beta },
\end{equation*}
hence by Gr\"{o}nwall inequality and the induction argument, it is
enough to control $\left\vert x\left( \cdot ,t\right) \right\vert
_{\ast }$ and $ \left\vert x_{\Gamma }\left( \cdot ,t\right)
\right\vert _{\beta }$, to guarantee that $x\left( \Gamma ,t\right)
\in C^{n,\beta}\left(
 \Gamma_{0},\Gamma _{1} \right) $, for all $n\geq 1$, (and consequently in $
C^{\infty }\left(
 \Gamma_{0},\Gamma _{1} \right)$, whenever $x_0\in C^{\infty }\left(
 \Gamma_{0},\Gamma _{1} \right)\cap \left\{ \left\vert x\right\vert _{\ast
}>0\right\}$).
\end{proof}

\section{\label{sec:concl}Conclusions}
The 2D Euler-$\alpha $ model \cite
{a_CFHOTW99_ChanPipe,a_HMR98,a_HMR98b,a_H02_pA,a_MS03,a_C01} is an
inviscid regularization of the Euler equations. In \cite{a_OS01} it
has been shown the existence of a unique global weak solution of 2D
Euler-$\alpha $ equations, when the initial vorticity is in the
space of Radon measures on ${\mathbb{R}}^{2}$. The
Birkhoff-Rott-$\alpha$ equation for the evolution of the 2D vortex
sheet is induced by the 2D Euler-$\alpha$ equations, and it is an
$\alpha$-analogue of the Birkhoff-Rott equation, induced by the 2D
Euler equations.

The structure of weak solutions of 2D Euler equations, for the
vortex sheet initial data with initial vorticity being a Radon
measure of a distinguished sign, given by Delort
\cite{a_D91,a_EM94,a_M93,a_LX95, a_S95, a_S96, b_MB02} is not known,
yet the BR equations assume \textit{a priori} that a vortex sheet
remains a curve at a later time. On the contrary, in the 2D
Euler-$\alpha$ case, if the vorticity is initially supported on a
curve, it remains supported on a curve for all times, hence the
BR-$\alpha$ equation gives an equivalent description of the vortex
sheet evolution, as the 2D Euler-$\alpha$ equations.

In this paper we report the global regularity of the BR-$\alpha$
approximation for the 2D vortex sheet evolution. We show that
initially smooth  self-avoiding vortex sheet remains smooth for all
times, under the condition that the initial density is an integrable
function of the vortex curve with respect to the arc-length measure.

Unlike the original BR problem that exhibits the Kelvin-Helmholtz
instability, the linearized, about the flat solution, BR-$\alpha$
model has growth rates that decay to zero for large wave numbers,
larger than $O(\alpha)$. This, in turn, is also an indication of the
role that the parameter $\alpha$ plays in slowing the process of
formation of scales smaller than $\alpha$. Another indication that
$\alpha$ controls the development of small scales, smaller than
$\alpha$, arises from the Lagrangian description of the flow. The
lower bound \eqref{eq:trajBelowBound} implies that the evolution of
small scales, relative to $\alpha$, at each instant of time, is
controlled from below by the initial ratio. That is, for any finite
time, the spatial scales smaller than alpha develop at a controlled
rate.

The linear stability analysis also implicates that the BR-$\alpha$
approximation could be closer to the original BR equation than the
existing regularizing methods, such as vortex blob model, due to the
less regular kernel. A numerical study comparing the $\alpha$ and
the vortex blob regularizations for planar and axisymmetric vortex
filaments and sheets is reported in \cite{a_HNP06}.

The full details of the results reported in this paper will be
presented in a forthcoming paper.

\begin{acknowledgments}
C.B. would like to thank the Faculty of Mathematics and Computer
Science at the Weizmann Institute of Science for the kind
hospitality where this work was initiated. This work was supported
in part by the BSF grant no.~2004271, the ISF grant no.~120/06, and
the NSF grants no.~DMS-0504619 and no.~DMS-0708832.
\end{acknowledgments}

\bibliography{2007_BLT_BRalpha}

% ------------------------------------------------------------------------
%Included for Gather Purpose only:
%input "2007_BLT_BRalpha.bib"

\end{document}